\documentclass[11pt]{article}
\usepackage{amssymb,amsmath,amsthm,amsfonts}

\newcommand{\gl}{\mathfrak{gl}}
\newcommand{\gee}{\mathfrak{g}}
\newcommand{\geven}{\mathfrak{g}_{\bar{0}}}
\newcommand{\godd}{\mathfrak{g}_{\bar{1}}}
\newcommand{\yoo}{\mathfrak{U}}
\newcommand{\pee}{\mathfrak{p}}
\newcommand{\hee}{\mathfrak{h}}
\newcommand{\kay}{\mathfrak{k}}

\newcommand{\zed}{\mathbb{Z}}
\newcommand{\zz}{\mathbb{Z}_2}
\newcommand{\zee}{\mathcal{Z}}
\newcommand{\C}{\mathbb{C}}

\newcommand{\emm}{\mathcal{M}}

\newcommand{\Veven}{V_{\bar{0}}}
\newcommand{\Vodd}{V_{\bar{1}}}
\newcommand{\en}{\mathfrak{n}}

\newtheorem{thm}{Theorem}[section]
\newtheorem{pro}[thm]{Proposition}

\newtheorem{cor}[thm]{Corollary}
  
\newtheorem{defn}[thm]{Definition}

\begin{document}
 \title{On Some Representations of Nilpotent Lie Algebras and Superalgebras} 
 \author{Shantala Mukherjee
 \thanks{This work is part of the author's doctoral dissertation, written
at the University of Wisconsin-Madison under the supervision of Prof.
Georgia Benkart, and financially  
supported in part by NSF grant \#{}DMS-0245082}\\
 Dept. of Mathematics \\
 DePaul University \\
 Chicago, IL 60614}

 \date{}
 \maketitle

 \begin{abstract}
 Let $G$ be a simply connected, nilpotent Lie group with Lie algebra $\gee$. The group $G$
acts on the dual space $\gee^*$ by the coadjoint action.
By the orbit method of Kirillov, the simple unitary
representations of
$G$ are in bijective correspondence
with the coadjoint orbits in $\gee^*$, which in turn are in bijective correspondence
with the primitive ideals of the universal enveloping algebra of
$\gee$. The
number of simple $\gee$-modules which have
a common eigenvector for a particular subalgebra of $\gee$ and are annihilated by
a particular primitive ideal $I$ is shown by Benoist to depend on geometric properties of
a certain subvariety of the coadjoint orbit corresponding to $I$.  
We determine the 
exact number of such modules  when the coadjoint orbit is two-dimensional.

Bell and Musson showed that the algebras obtained by factoring the
universal enveloping superalgebra of a Lie superalgebra by
graded-primitive ideals are isomorphic to
tensor products of Weyl algebras and Clifford algebras. We describe certain cases where
the factors are purely Weyl algebras and determine how the sizes of these Weyl algebras depend
on the graded-primitive ideals.

 \end{abstract}

 \section{Introduction}
 The study of the representations of a real Lie group $G$ is related to
the study of the representations of its complex Lie algebra $\gee$.
If $G$ is simply
connected and nilpotent, then the irreducible unitary representations of
$G$ are related to certain ideals of the universal enveloping algebra
of $\gee$. This correspondence links ideas from geometry, analysis
and algebra. The representation theory of Lie superalgebras is
similar to, yet different from the representation theory of Lie
algebras. Lie superalgebras are of interest to physicists in the
context of supergauge symmetries relating particles of different statistics.
This paper is devoted to the study of some aspects of the representation theory of
nilpotent Lie algebras and superalgebras.

The orbit method was created by Kirillov in the attempt to describe the unitary dual $\hat{N_m}$ for the nilpotent Lie group $N_m$ of
$m \times m$ upper triangular matrices with $1$'s on the diagonal (the unitriangular group).
It turned out that the orbit method had much wider applications. In Kirillov's words:
`\dots all main questions of representation theory of Lie groups: construction of irreducible
representations, restriction-induction functors, generalized and infinitesimal characters, Plancherel
measure, etc., admit a transparent description in terms of coadjoint orbits' ({\cite{Kirillov:orb03}}).

For a nilpotent Lie algebra $\gee$,
Dixmier in {\cite{Dixmier:ea96}} formulated the correspondence between the set of primitive
ideals of the universal enveloping algebra $U(\gee)$, the set of coadjoint orbits in
$\gee^*$, and the sizes of the Weyl algebras obtained by factoring $U(\gee)$ by primitive ideals.
In general, for any primitive ideal $I$, there are infinitely many non-isomorphic simple $\gee$-modules that have $I$ as their annihilator.
Benoist in {\cite{Benoist:ens90}} used the orbit method to show, for a finite-dimensional nilpotent
Lie algebra $\gee$, that the  number (of isomorphism classes) of simple $\gee$-modules which are annihilated by
a primitive ideal $I$ \emph{and} have a common eigenvector for a certain subalgebra of $\gee$, depends
on properties of a certain subvariety of the coadjoint orbit corresponding to $I$ under the Dixmier
correspondence. We refine his results in a particular case.

The methods used by Dixmier and Conze ({\cite{Conze:cr71}}) to describe the primitive ideals of
universal enveloping algebras have been extended to
study graded-primitive ideals of universal enveloping superalgebras of Lie superalgebras  by Letzter
({\cite{Letzter:cmams92}}), and by Bell and Musson ({\cite{BellMusson:lms90}}, {\cite{Musson:am92}}). 
In {\cite{BellMusson:lms90}} it is shown that the algebras obtained by factoring the universal enveloping superalgebra by
graded-primitive ideals are isomorphic to
tensor products of Weyl algebras and Clifford algebras. We describe certain cases where
the factors are purely Weyl algebras and determine how the sizes of these Weyl algebras depend
on the graded-primitive ideals.

\section{ Basic Definitions for Lie Algebras }
Our objects of study in sections 2 through 9 are finite-dimensional nilpotent Lie algebras $\gee$ over the field $\C$
of complex numbers.
Thus if
\[ \gee^0 = \gee,\ \gee^1 = [\gee,\gee],\ \gee^2 = [\gee,\gee^1],\dots,\ \gee^{i+1} = [\gee,\gee^i] .\]
\noindent is the lower central series of $\gee$, then $\gee^k = 0$ for some $k \geq 1$.

By the Birkhoff embedding theorem
({\cite[Thm.~1.1.11]{Corwin:nlg90}}), any nilpotent Lie algebra of
finite dimension
over $\C$ is isomorphic to a subalgebra of $\en_m$ for some $m$, where $\en_m$ is the Lie algebra of
strictly upper triangular $m \times m$ matrices under the product $[x,y]=xy-yx$ for $x,y \in \en_m$.

Let $G$ be a connected, simply connected nilpotent Lie group with Lie algebra $\gee$. If $\exp$
is the exponential map from $\gee$ to $G$, then $\exp(\gee) = G$ ({\cite[Thm.~1.2.1]{Corwin:nlg90}}).
In fact, if we identify $\gee$ with a subalgebra of $\en_m$, then the exponential map becomes the ordinary exponential map $x \mapsto
\sum_{j=0}^{\infty}\frac{1}{j!}x^j$.
The group $G$ acts on
$\gee$ by the \emph{adjoint} action,
\[ \textrm{Ad}g(x) = gxg^{-1},\ \ \forall x \in \gee,\ g \in G, \]
\noindent and it
acts on the dual space $\mathfrak{g}^*$
by the \emph{coadjoint} action. Thus, if $f \in \gee^*$, then
\[ (g.f)(y) = f(g^{-1}yg),\ \ \ \forall g \in G, y \in \gee. \]

The orbit of $f \in \gee^*$ under the action of $G$ is called the
 \emph{coadjoint orbit}
containing $f$. It is denoted by $\Omega_f$.

\section{ Primitive Ideals and Weyl Algebras}

An ideal $I$ of  the universal enveloping algebra $U=U(\gee)$ is said to be \emph{primitive} if it is the annihilator of a simple
left $\gee$-module. The set of primitive ideals of $U(\gee)$ is denoted by $\textrm{Prim}~U$.
For a nilpotent Lie algebra $\gee$, any primitive ideal of $U(\gee)$ is maximal among
the set of proper two-sided ideals of $U(\gee)$ by {\cite[Prop. 4.7.4]{Dixmier:ea96}}.

\begin{defn}\label{D:weylalg}
For $n\geq 1$, the $n$-th \emph{Weyl algebra} $\mathcal{A}_n$ is the algebra with $2n$ generators
$p_1, q_1,\dots, p_n, q_n$, and relations
\[ [p_i, q_i ] = 1, \]
\[ [p_i, q_j] = [p_i, p_j] = [q_i, q_j] = 0\ \  \text{for}\ i \neq
j. \]
By convention $\mathcal{A}_0 = \C$.
\end{defn}

By {\cite[Prop. 4.7.9]{Dixmier:ea96}}, if $I$ is a primitive ideal of $U(\gee)$, then
the quotient $U(\gee) / I$ is isomorphic to $\mathcal{A}_n$, for some positive integer $n$.
This integer $n$ is uniquely determined by the ideal $I$, and it is called the
\emph{weight} of $I$ ({\cite[4.7.10]{Dixmier:ea96}}).

\section{Coadjoint Orbits and Primitive Ideals}\label{S:orbitideal}

Any $f \in \gee^*$ determines an alternating bilinear form $B_f$ on $\gee$ given by
\[ (x,y) \mapsto B_f(x,y) := f([x,y])\ \ \ \forall x,y \in \gee. \]
Let
\begin{align}
\gee^f = &\{x \in \gee \mid f([x,y])=0 \ \forall y \in \gee\}  \notag \\
           = &\{x \in \gee \mid B_f(x,y) = 0\ \forall y \in \gee \}. \notag
\end{align}
It is obvious that $\gee^f$ is a subalgebra of $\gee$ by the Jacobi identity on $\gee$. It is called the \emph{radical} of $f$,
or the \emph{kernel} of the form $B_f$.

A Lie subalgebra $\kay$ is said to be \emph{subordinate} to $f$ if $f([x,y])=0$ for all $x,y \in \kay$,
i.e. if $\kay$ is a totally isotropic subspace of $\gee$ with respect to the alternating bilinear form
$B_f$. The largest dimension of a subalgebra subordinate to $f$ is
$\frac12(\dim\gee + \dim\gee^f)$ by {\cite[1.12.1]{Dixmier:ea96}}. A subalgebra that is subordinate
to $f$  and has this maximal dimension is called a \emph{polarisation} of $f$.
The set of all polarisations of $f$ is denoted by $P(f)$. If $\pee \in P(f)$, then $\pee \supseteq \gee^f$
(see {\cite[1.12.1]{Dixmier:ea96}}).

Let $\hee$ be a Lie subalgebra of $\gee$ and let $W$ be an $\hee$-module.  Since
$U(\gee)$ is a right $U(\hee)$-module under multiplication, we can form the induced module
$V = U(\gee)\otimes_{U(\hee)}W$  with $U(\gee)$-action given by
\[ x(u\otimes w) = xu\otimes w\ \ \ \forall x,u \in U(\gee)\ \text{and}\ w \in W.\]

Assume $f \in \gee^*$ and $\pee \in P(f)$.  Then, $f([x,y]) = 0$ for all $x,y \in \pee$, so we can define
a $\pee$-module action on a one-dimensional vector space $\{\pee,f\}=\C v$ as follows:
\[ x.v = f(x)v,\ \ \forall x \in \pee. \]
Then $Ind_{\pee}^{\gee}\{\pee,f\}=U(\gee)\otimes_{U(\pee)}\C v$, the
$\gee$-module induced from $\{\pee,f\}$, is a simple
$\gee$-module, hence a simple $U(\gee)$-module by
{\cite[Thm. 6.1.1]{Dixmier:ea96}}.
The annihilator in $U(\gee)$ of this module is a primitive ideal, denoted by $I(f)$.
By {\cite[Thm. 6.1.4]{Dixmier:ea96}}, the ideal $I(f)$ depends only on $f$, and not on the choice of
the polarisation $\pee \in P(f)$. Any primitive ideal $I$ is of the
form $I(f)$, for some $f \in \gee^*$
({\cite[Sec.~6.1.5, Thm.~6.1.7]{Dixmier:ea96}}).

If $f \in \gee^*$ and $a \in G$, then $I(a.f) = I(f)$ by {\cite[Prop. 2.4.17]{Dixmier:ea96}}.
Thus, the map
\[ f \mapsto I(f), \]
between $\gee^*$ and $\textrm{Prim}~U$, defines a map
\[ \Omega_f  \mapsto I(f) \]
 between the set of $G$-coadjoint orbits
in $\gee^*$ and $\textrm{Prim}~U$. By {\cite[Thm. 6.1.7, Prop. 6.2.3, Thm. 6.2.4]{Dixmier:ea96}},
this map is a bijection. Moreover, we have
$U(\gee) / I(f) \simeq \mathcal{A}_n$, where
$n = \frac{1}{2}\text{rank}(B_f )=\frac{1}{2}\dim(\gee/\gee^f)$ according to {\cite[Prop. 6.2.2]{Dixmier:ea96}}.
Notice that if $\gee$ is a finite-dimensional nilpotent Lie algebra,
and if $f \in \gee^*$, then the weight of the associated primitive
ideal $I(f)$ equals $\dim (\gee / \pee)$ for any polarisation $\pee$ of
$f$.

\section{Geometry of Coadjoint Orbits}\label{S:orbgeom}
First we recall some generalities about symplectic vector spaces and symplectic varieties,
and then we describe the symplectic structure on a coadjoint orbit.

A symplectic structure on an even-dimensional vector space $V$ is determined by a non-degenerate
alternating bilinear form $\omega$ on $V$. If $W$ is a subspace of $V$, then its
\emph{orthogonal complement} is
\[ W^{\bot} = \{v \in V \mid \omega(w,v)=0\ \forall w \in W \}. \]
If $W \subseteq W^{\bot}$, then $W$ is called \emph{isotropic}.
If $W^{\bot} \subseteq W$, then $W$ is called \emph{coisotropic}.
A subspace $W$ that is both isotropic and coisotropic is said to be a \emph{lagrangian} subspace.
A lagrangian subspace of $V$ is always of dimension $\frac{1}{2}\dim V$ (see {\cite[Ch.~11, Prop.~2.2]{Coutinho:dmod95}}).

Let $\mathcal{V}$ be an algebraic variety. A symplectic structure on $\mathcal{V}$ is a
non-degenerate algebraic 2-form $\omega$ on $\mathcal{V}$ such that $d\omega = 0$
(see {\cite[Sec.~1.1]{Chriss:comgeo97}} for details). If $p \in \mathcal{V}$, then there is a alternating bilinear form $\omega_{p}$
on $T_{p}\mathcal{V}$, the tangent space at $p$.
If $\mathcal{W}$ is a subvariety of $\mathcal{V}$, then it is said to be \emph{lagrangian} if the tangent
space $T_p\mathcal{W}$ is a lagrangian subspace of $T_p\mathcal{V}$ at every non-singular point $p \in \mathcal{V}$.
The dimension $\dim\mathcal{V}$ of a variety $\mathcal{V}$ is defined to be the dimension of
the tangent space $T_p\mathcal{V}$ at any non-singular point $p \in \mathcal{V}$. Thus,
if $\mathcal{W}$ is a lagrangian subvariety of $\mathcal{V}$, then $\dim\mathcal{W} = \frac{1}{2}\dim\mathcal{V}$.

Any coadjoint orbit $\Omega \subset \gee^*$ has a natural symplectic
structure given as follows (see {\cite[Prop.~1.1.5]{Chriss:comgeo97}}):

Assume $f \in \Omega$. The tangent space $T_f\Omega$ at the point $f$
is equal to $\gee / \gee^f$. We define an alternating bilinear form on $\gee$
\[ \omega_f \colon \gee \times \gee \rightarrow \C,\ \ \omega_f \colon
(x,y) \mapsto B_f(x,y) = f([x,y]).\]
The form $\omega_f$ descends to $\gee/\gee^f$. Thus the assignment
$f \mapsto \omega_f$ gives a non-degenerate 2-form $\omega$ on $\Omega$
such that $d\omega = 0$.

Notice that $\dim\Omega = \dim(\gee/\gee^f) = \text{rank}(B_f)$ which is an even number.
Any coadjoint orbit is an irreducible variety of $\gee^*$ (see {\cite[Prop.~8.2]{Hum:linalg75}}).

\section{ Generalized Weight Modules }

In this section, we relax our assumptions and let
$\mathfrak{g}$ be an arbitrary finite-dimensional Lie algebra and
$\mathfrak{h}$ be a nilpotent Lie subalgebra of $\mathfrak{g}$.

\begin{defn}
A $\mathfrak{g}$-module $N$ is a \emph{ generalized weight module over}
$\mathfrak{h}$, if as an $\mathfrak{h}$-module, $N$ decomposes as
\[ N = \bigoplus_{\mu \in suppN}N^{\mu} \]
for some subset supp$N$ of $\mathfrak{h}^*$, where $N^{\mu}$ is
a non-zero generalized weight space of weight $\mu \in \mathfrak{h}^*$ for each $\mu \in \textrm{supp}N$, i.e. for
every $v \in N$ and every $x \in \hee$, there exists $l=l(x,v) \in \zed_{>0})$ such that
$(x - \mu(x))^{l}.v = 0$.
\end{defn}
As an $\hee$-module under the adjoint action, $\gee$ decomposes as
\[ \gee = \hee' \oplus \left (\bigoplus_{\alpha \in \Delta_{\hee}\subset \hee^*}\gee^{\alpha}\right ), \]
\noindent where each $\gee^{\alpha}$ is a generalized weight space of weight $\alpha \in
\hee^*$, $\alpha \neq 0$ (i.e., for every $x \in \gee^{\alpha}$ and every $y \in \hee$ there exists
$n \in \zed_{\geq 0}$ so that $(ad(y) - \alpha(y))^nx = 0$), and $\hee'$ is the generalized weight
space of weight $0$. The set $\Delta_{\hee}$ of all non-zero weights $\alpha$ is the set of
$\hee$-roots of $\gee$.

An element $x \in \mathfrak{g}$ is said to act \emph{ locally finitely}
on a $\mathfrak{g}$-module $N$ if the vector space spanned by the
vectors $v, x.v, x^2.v,\dots $ is finite-dimensional for any $v \in N$.
By {\cite[Cor.~2.7]{Fernando:tams90}},
the set $\gee[N]$ of all elements  of $\mathfrak{g}$ which act locally finitely on $N$
is a Lie subalgebra  of $\mathfrak{g}$. It is the largest subalgebra
of $\gee$ that is locally finite on $N$ and is called
the \emph{ Fernando subalgebra } of $\mathfrak{g}$ with respect to $N$. For any $v \in N$,
the $\gee[N]$-module
$U(\mathfrak{g}[N]).v$ generated by $v$ is finite-dimensional by {\cite[Prop.~1]{PenkovSerganova:ja98}}.

Let $N$ be a $\gee$-module, possibly infinite-dimensional. We say that $x \in \gee$ acts
\emph{freely} on $N$ if the vectors $v, x.v, x^2.v,\dots$ are linearly
independent for any $v \in N$. If $N$ is a simple
 $\mathfrak{g}$-module, then any $x \in \mathfrak{g}$ acts either locally finitely or   
freely on $N$.

For a nilpotent Lie subalgebra $\mathfrak{h} \subset \mathfrak{g}$, $N$ is
a generalized weight module over $\mathfrak{h}$ precisely when
$\mathfrak{h} \subset \mathfrak{g}[N]$.
If $N$ is a simple $\gee$-module, then we denote by $\Gamma_N$ the cone in
$\langle \Delta_{\hee} \rangle _{\mathbb{R}}$, the real span of the
roots, generated by all $\alpha \in \Delta_{\hee}$ such that
$\gee^{\alpha}$ is not contained in $\gee[N]$.

If $\gee$ is a finite-dimensional solvable Lie algebra and $N$ is a
simple $\gee$-module, then $\textrm{supp}N = \mu + \Gamma_N$, for some
$\mu \in \hee^*$ ({\cite[Rem., Prop.~2]{PenkovSerganova:ja98}}).
If, in addition, $\gee$ is nilpotent, then the set $\Delta_{\hee}$ is empty by Engel's theorem; therefore
$\Gamma_N = 0$.  To summarise: if $\gee$ is a nilpotent Lie algebra, $\hee$ is a
subalgebra of $\gee$, and $N$ is a simple $\gee$-module that is a generalized weight module
over $\hee$, then $\textrm{supp}N = \{\mu \}$, for some $\mu \in \hee^*$.

\section{ Simple Modules Containing an Eigenvector \\ for a Subalgebra }\label{S:sympvar}

Let $\gee$ be a nilpotent Lie algebra and let $f \in \gee^*$. Assume $\kay$  is a subalgebra
of $\gee$ subordinate to $f$, and let
$N$ be a $\gee$-module. Thus $f([\kay,\kay]) = 0$.
Set
\begin{equation} \kay(f) = \{x - f(x) \mid x \in \kay\} \subset U(\gee), \end{equation}
\begin{equation} N^{\kay, f} = \{ n \in N \mid (x - f(x)).n = 0 \ \
\forall x \in \kay \}, \end{equation}
and
\begin{equation} \kay^{\top} = \{\lambda \in \gee^* \mid \lambda(\kay) = 0\}. \end{equation}
The set $f + \kay^{\top} = \{\lambda \in \gee^* \mid \lambda(x)=f(x)\ \forall x \in \kay\}$ is an affine linear subspace of $\gee^*$. We will show that it is an irreducible algebraic variety.
Let $\dim \gee = r$ and $\dim \kay = s,\ s\leq r$. Suppose $e_1,e_2,\dots,e_r$ is a basis of $\gee$ such that $e_1,e_2,\dots,e_s$ is a basis of $\kay$.
 We have
the dual basis $e_1^*,e_2^*,\dots,e_r^*$ of $\gee^*$. Any $\lambda \in \gee^*$ can be uniquely
represented by $(a_1,\dots,a_r) \in \C^r$ where $\lambda = \sum_{1\leq i\leq r}a_ie_i^*$.
Thus
\begin{align}
f+\kay^{\top} &= \{\lambda \in \gee^*\mid \lambda(e_i)=f(e_i),\ i=1,\dots,s\} \notag \\
&=\{(a_1,\dots,a_s,\dots,a_r)\mid a_i=f(e_i),\ i=1,\dots,s\}. \notag
\end{align}
So the coordinate ring of $f+\kay^{\top}$ is $\C[a_1,\dots,a_r]/\langle a_1-f(e_1),\dots,a_s-f(e_s)\rangle
\simeq \C[\overline{a_{s+1}},\dots,\overline{a_r}]$ which does not have zero divisors, hence
$f+\kay^{\top}$ is irreducible (see {\cite[Prop.~2.2.5]{Chriss:comgeo97}}).
Notice that $\dim(f+\kay^{\top}) = \dim \kay^{\top} = \dim(\gee/\kay)$.

We will use the following theorem and corollary of Benoist ({\cite[Thm.~6.1]{Benoist:ens90}}). Recall
the definitions of lagrangian spaces and varieties from Section~\ref{S:orbgeom}.
\begin{thm}\label{BenLag}
Let $\mathfrak{g}$ be a finite-dimensional nilpotent Lie algebra over $\mathbb{C}$, and let
$U = U(\mathfrak{g})$. Assume $\Omega $ is a $G$-orbit in $\mathfrak{g}^*$, and $I$ is the
primitive ideal of $U$ associated with $\Omega$. Let $f \in \mathfrak{g}^*$ and
$\mathfrak{k}$ be a Lie subalgebra of $\gee$ such that $f([\mathfrak{k},\mathfrak{k}])
= 0$. Set $\zee = \Omega \cap (f + \mathfrak{k}^{\top})$ and let $M$ be the
$\mathfrak{g}$-module $U/(I + U\mathfrak{k}(f))$. Assume $S$ is the set of
those simple $\mathfrak{g}$-modules $N$ with annihilator $I$ such
that $N^{\mathfrak{k},f} \neq 0$.
\begin{enumerate}
\item
The following are equivalent:
\begin{description}
\item[(i)]
$\zee$ is a lagrangian subvariety of $\Omega$.
\item[(ii)]
$M$ is of finite length (has a finite composition series).
\item[(iii)]
$S$ is a finite set.
\end{description}
\item
If one (hence all) of the conditions in part 1. hold, then:
\begin{description}
\item[(a)]
$\zee$ is a smooth variety.
\item[(b)]
There is a bijection  between
the irreducible components $\Lambda$ of $\zee$ and the elements $M_{\Lambda}$ of $S$.
\item[(c)]
There is an
isomorphism of $\mathfrak{g}$-modules
$M \simeq \bigoplus_{\Lambda}M_{\Lambda}^{\oplus m_{\Lambda}}$, where
$m_{\Lambda} = \dim(M_{\Lambda}^{\kay, f})$. In particular,
$m_{\Lambda}$ is finite for each $\Lambda$ and $M$ is semi-simple.
\end{description}
\end{enumerate}
\end{thm}

\begin{cor}\label{BenLagcor}
With assumptions as in Theorem~\ref{BenLag},  
\begin{enumerate}
\item
$\zee = \emptyset \Longleftrightarrow M = 0 \Longleftrightarrow S = \emptyset. $
\item
If $\zee$ is an orbit under the group $K = \exp(\mathfrak{k})$, then $M$ is a multiple
of the simple module $M_{\zee}$, $M = M_{\zee}^{\oplus m_{\zee}}$.
\end{enumerate}
\end{cor}

Since $\mathfrak{g}$ is assumed to be a nilpotent Lie algebra, all of
its
subalgebras are nilpotent too. Let $N$ be a member of the set $S$ as
defined in Theorem~\ref{BenLag}. Then the subalgebra $\mathfrak{k}$ acts
locally finitely on some element of $N$, hence on all of $N$. Thus
$\mathfrak{k} \subset \mathfrak{g}[N]$. So $N$ is a generalized weight
module over $\mathfrak{k}$ and supp$N$ = $f\vert _{\mathfrak{k}}$.

On the other hand, suppose $N$ is a simple generalized weight module
over some subalgebra $\mathfrak{k} \subset \mathfrak{g}$ with
supp$N = \{\mu\}$ relative to $\kay$. Let $I =\textrm{ann}_{U(\gee)}(N)$, and let $\mu = f\vert _{\mathfrak{k}}$
for some $f \in \mathfrak{g}^*$. Let $v \in N$ be any non-zero generalized
weight vector of weight $\mu$ relative to $\kay$. Then
$U(\mathfrak{g}[N])v$ is a finite-dimensional $\mathfrak{g}[N]$-module.
Let $N_v$ be a simple $\mathfrak{g}[N]$-submodule of $U(\gee[N])v$. Since
$\mathfrak{g}[N] \subseteq \gee$ is nilpotent,
$N_v$ is a one-dimensional $\gee[N]$-module, hence a one-dimensional $\mathfrak{k}$-module, and $x \in \kay$ acts by multiplication
by $f(x)$ on it, which implies that
$f([\mathfrak{k},\mathfrak{k}]) = 0$ and also that
$N^{\mathfrak{k},f} \neq 0$.
Thus, we have proved the following result:

\begin{thm}\label{genwt}
Assume $\mathfrak{g}$ is a finite-dimensional nilpotent
 Lie algebra over $\mathbb{C}$,
and let $U = U(\gee)$. Let $I \in \textrm{Prim}~U, f \in
\mathfrak{g}^* $, and $\mathfrak{k}$ be a Lie subalgebra of $\mathfrak{g}$
such that $f([\mathfrak{k},\mathfrak{k}]) = 0$. Let
$\mathcal{E}_I^{\mathfrak{k},f}$  be the set of all simple
$\mathfrak{g}$-modules $N$ with annihilator $I$ such that $N^{\mathfrak{k},f} \neq 0$,
and let $\mathcal{W}_I^{\mathfrak{k},f}$  be the set of
all simple $\mathfrak{g}$-modules $N$ with annihilator I which are
also generalized weight modules over $\mathfrak{k}$ with supp$N =
\{ f\vert_{\mathfrak{k}} \}$. Then
\[ \mathcal{E}_I^{\mathfrak{k},f} = \mathcal{W}_I^{\mathfrak{k},f}. \]
\end{thm} \qed

As a consequence we have
\begin{cor}
With assumptions as in Theorem~\ref{genwt}, let $\kay = \pee$, a
polarisation of $f$ in $\gee$. If $\Omega$ is a $G$-orbit in $\gee^*$
and $f \notin \Omega$, then $\mathcal{E}_I^{\pee,f}=\emptyset$, hence $\mathcal{W}_I^{\pee,f} = \emptyset$.
\end{cor}

\noindent \textbf{Proof.}
If $\Omega$ is a $G$-orbit in
$\gee^*$ and $f \notin \Omega$, then $f + \pee^{\top}$ is contained in the coadjoint orbit of
$f$ by {\cite[Prop.~1, Sec.~2]{Kirillov:bams99}}. Since distinct coadjoint orbits are disjoint,
 $\zee = \Omega \cap (f + \mathfrak{p}^{\top}) = \emptyset$, which,
from Corollary~\ref{BenLagcor} and Theorem~\ref{genwt} implies that
$\mathcal{E}_I^{\pee,f}=\mathcal{W}_I^{\mathfrak{p},f} = \emptyset$. \qed
\bigskip

Next we consider what happens when $\Omega
= \Omega_f$, the coadjoint orbit passing through $f$ itself, and
$\mathfrak{k} = \mathfrak{p} \in P(f)$, a polarisation of $f$. Then $\zee = \Omega_f \cap
(f + \mathfrak{p}^{\top}) = f + \mathfrak{p}^{\top}$, by
{\cite[Prop.~1, Sec.~2]{Kirillov:bams99}}. We want to determine how many irreducible components
the variety $\zee$ has when $\zee$ is lagrangian.

\section{ Induced Modules }
Assume $\gee$ is a nilpotent Lie algebra.
Let $f \in \mathfrak{g}^*$ be such that $f([\mathfrak{g}, \mathfrak{g}])
\neq 0$. Let $\mathfrak{p} \in P(f)$ be a polarisation of $f$. Let $\{\pee,f\}$ be the one-dimensional
$\mathfrak{p}$-module $\C v$ given by $f \vert_\mathfrak{p}$.
By {\cite[Prop.~6.2.9]{Dixmier:ea96}, the induced
$\mathfrak{g}$-module $\mathcal{M}_f := Ind_{\pee}^{\gee}\{\pee,f\}=U(\mathfrak{g})\otimes_{U(\mathfrak{p})}\{\pee,f\}$
is simple. The mapping $v \mapsto 1\otimes v$ of $\{\pee,f\}$ into $\emm_f$ is an
injective $\mathfrak{p}$-module homomorphism. So $\{\pee,f\}$ can
be identified with a submodule of the $\pee$-module $\emm_f$ under this mapping.

As mentioned in Section 2.3, the primitive ideal $I$ associated with the
coadjoint
orbit passing through $f$ is the precisely the annihilator of the
module $\emm_f$ in $U$, and it does not depend on the choice of polarisation
of $f$.

Let $J = \textrm{ann}_{U(\mathfrak{p})}(\{\pee,f\})$.
From {\cite[Prop.~5.1.7]{Dixmier:ea96}}, we have
\[ \textrm{ann}_{U(\mathfrak{g})}(\{\pee,f\}) = U(\mathfrak{g})J, \]
which is a left ideal of $U(\mathfrak{g})$. And,
 $\textrm{ann}_{U(\mathfrak{g})}(\emm_f) = I $ is the largest two-sided ideal of
$U(\mathfrak{g})$ contained in $U(\mathfrak{g})J$.

By Prop.~5.1.9 (i) in {\cite{Dixmier:ea96}},
the mapping $\phi$ of $U(\mathfrak{g})$
into $\emm_f$ defined by $\phi(u) = u \otimes v$ for all $u \in
U(\mathfrak{g})$
is surjective and has kernel $U(\mathfrak{g})J$.

By Prop.~5.1.9 (ii) in {\cite{Dixmier:ea96}}, the mapping $\bar{\phi}$ of
$U(\mathfrak{g})/U(\mathfrak{g})J$ into $\emm_f$ inherited from $\phi$ by
passage to the quotient is a $\mathfrak{g}$-module isomorphism.
Thus $U(\mathfrak{g})/U(\mathfrak{g})J \simeq \emm_f$ as
$\mathfrak{g}$-modules, and so
by Prop.~5.1.9 (iii) in {\cite{Dixmier:ea96}}, we
see that $U(\mathfrak{g})J = U(\mathfrak{g}){\mathfrak{p}}(f)$,
where
\[\pee(f) = \{x - f(x) \mid x \in \pee\} \subseteq U(\gee).\]
Consequently, $U(\mathfrak{g})/U(\mathfrak{g}){\mathfrak{p}}(f) \simeq \emm_f$
and so is a simple $\mathfrak{g}$-module, since $\emm_f$ is a simple $\gee$-module.
But from Prop.~5.1.7 (ii) in {\cite{Dixmier:ea96}}, we have
$I \subset U(\mathfrak{g}){\mathfrak{p}}(f)$. Hence $I+U(\gee)\pee(f) = U(\gee)\pee(f)$ and
thus, $U(\mathfrak{g})/(I + U(\mathfrak{g}){\mathfrak{p}}(f))=U(\gee)/U(\gee)\pee(f)$ is
a simple $\mathfrak{g}$-module.

As shown in Section~\ref{S:sympvar}, $\zee = \Omega_f \cap (f + \pee^{\top})=f+\pee^{\top}$ is irreducible, hence it has
only one irreducible component.
Therefore if $\zee$ is a lagrangian subvariety of $\Omega_f$, then the set $\mathcal{E}_I^{\mathfrak{p},f}$
contains only one element $M_{\zee}$. We see that
$U(\mathfrak{g})\otimes_{U(\mathfrak{p})}{\{\pee,f\}}$ is an element of
$\mathcal{E}_I^{\mathfrak{p},f}$, so in this case it is isomorphic to
$M_{\zee}$.

\begin{thm}\label{twodimlag}
If $f \in \mathfrak{g}^*$ and the coadjoint orbit $\Omega_f$ passing
through $f$ is two-dimensional, then for the primitive ideal $I$ associated with $\Omega_f$,
the set $\mathcal{E}_I^{\mathfrak{p},f}$ contains only one element, $U(\gee)\otimes_{U(\pee)}\{\pee,f\}$,
for any $\pee \in P(f)$.
\end{thm}

\noindent \textbf{Proof.}
If $\dim \Omega_f = 2$, then $\dim{\gee} - \dim \mathfrak{p} = 1$, so the irreducible variety
$\zee = f + \mathfrak{p}^{\top}$ is one-dimensional, hence lagrangian. Thus the set
$\mathcal{E}_I^{\pee,f}$ contains only one element $U(\gee)\otimes_{U(\pee)}\{\pee,f\}$.   
\qed

\section{ Example: $\mathfrak{g} = \en_3$ }

Let $\mathfrak{g} = \en_3$, the Lie algebra of all $3 \times
3$ strictly
upper triangular matrices over $\mathbb{C}$. We can select a basis   
$\{x,y, z\}$ for $\gee$, where $x = E_{12}, y = E_{23}, z = E_{13}$ are standard matrix units.
Then, $[x,y] = z$, and $\C z$ is the center of $\gee$, so that $\en_3$ is isomorphic to
the Heisenberg Lie algebra.
Since
$[[\mathfrak{g},\mathfrak{g}],\mathfrak{g}] = 0$, for any $f
\in \mathfrak{g}^*$ there are only two possibilities:
\begin{enumerate}
\item
$f([\mathfrak{g},\mathfrak{g}]) = 0$, or
\item
$ f([\mathfrak{g},\mathfrak{g}])
\neq 0 $ and $f([[\mathfrak{g},\mathfrak{g}],\mathfrak{g}]) = 0$.
\end{enumerate}

In Case 1, $f$ defines a one-dimensional $\mathfrak{g}$-module
$\mathbb{C}v$ with $x.v = f(x)v$ for all $x \in \mathfrak{g}$.

In Case 2, it is easy to see that any proper subalgebra of $\mathfrak{g}$
is subordinate to $f$, and that any two-dimensional subalgebra of
$\mathfrak{g}$ is a polarisation of $f$. This means that the coadjoint
orbit of $f$ is two-dimensional (because the codimension of the polarisation is 1).
Let $\Omega$ be a coadjoint orbit in $\gee^*$, and
let $I \in \textrm{Prim}~U$ be the corresponding primitive ideal.
Now we divide considerations according to whether $f$ belongs to $\Omega$ or not.
\begin{description}
\item[(i)]
$f \in \Omega$
\end{description}
Let $\mathfrak{k}$ be a subalgebra of $\mathfrak{g}$ and suppose $f([\kay,\kay])=0$.
Recall that $\kay^{\top} = \{p \in \gee^* \mid p(\kay) = 0\} \subseteq \gee^*$.
Then $\zee = (f + \mathfrak{k}^{\top}) \cap \Omega \neq \emptyset$,
because $f \in \Omega$ and $f \in f+\kay^{\top}$.

If $\mathfrak{k} = 0$, then $\mathfrak{k}^{\top} = \mathfrak{g}^*$
and $\zee = \Omega$, so $\zee$ is not lagrangian.

If $\mathfrak{k} = \mathbb{C}z$, then $f + \mathfrak{k}^{\top}
= \Omega$, so $\zee = \Omega$ is not lagrangian.

If $\mathfrak{k} = \mathbb{C}(ax + by+ cz)$, where not both   
$a, b$ are zero, then $\zee = (f + \mathfrak{k}^{\bot}) \cap \Omega$
is a one-dimensional subvariety of the two-dimensional variety
$\Omega$, and so is lagrangian and irreducible. In this case,
the set $\mathcal{E}_I^{\mathfrak{k},f}$ has a unique element.   


If $\mathfrak{k}$ is two-dimensional, it is a polarisation of
$f$. Thus $\kay = \mathfrak{p} \in P(f)$, and then $f + \mathfrak{p}^{\top} \subset \Omega$.
Therefore $\zee = f + \mathfrak{p}^{\top}$ is lagrangian and irreducible,
so there is a unique
element $M_{\zee}$ in the set $\mathcal{E}_I^{\mathfrak{p},f}$ by Theorem~\ref{twodimlag} .
 Let $\{\pee,f\}$
be the one-dimensional $\mathfrak{p}$-module given by $f$. Then $M_{\zee}$ is
isomorphic to
the induced  irreducible $\mathfrak{g}$-module
$U(\mathfrak{g})\otimes_{U(\mathfrak{p})}\{\pee,f\}$, which we have shown
is isomorphic to the simple $\mathfrak{g}$-module
$M = U/(I + U\mathfrak{p}(f))$. By part 2 of Corollary 3.2,
we have $M$ isomorphic to ${M_{\zee}}^{\oplus m_{\zee}}$. But in this case $M$ is a
simple $\mathfrak{g}$-module, so it must be that $m_{\zee} = 1$.

\begin{description}
\item[(ii)] $f \notin \Omega$
\end{description}
Let $\mathfrak{k}$ be a subalgebra of $\mathfrak{g}$ such that
$f([\mathfrak{k},\mathfrak{k}]) = 0$.

If $\mathfrak{k} = 0$ then $\mathfrak{k}^{\top} = \mathfrak{g}^*,$
and so $\zee = \Omega$, which means $\zee$ is not lagrangian.

If $\mathfrak{k} = \mathbb{C}z$, then $\zee = \emptyset$, so $\mathcal{E}^{\kay,f}_I = \emptyset$.

If $\mathfrak{k} = \mathbb{C}(ax + by + cz)$ with not both
$a, b$ zero, then $\zee = (f + \mathfrak{k}^{\top}) \cap
\Omega$ is one-dimensional, hence lagrangian and irreducible.
So $\mathcal{E}_I^{\mathfrak{k},f}$ has a unique element.

If $\mathfrak{k} = \mathfrak{p}$, a polarisation of $f$, then
$(f + \mathfrak{p}^{\top}) \cap \Omega = \emptyset$, so
$\zee$ is empty, hence $\mathcal{E}^{\pee,f}_I = \emptyset$.

\section{Lie Superalgebras, Induced Modules and Graded-Primitive Ideals}
In this section, we assume that $\gee$ is a finite-dimensional Lie superalgebra over $\C$.
Thus $\gee$ has a $\zz$-grading, $\gee = \geven \oplus \godd$, and
a bilinear product $[ , ] : \gee \times \gee \rightarrow \gee$ such that
\begin{description}
\item[(1)] $[\gee_{\alpha}, \gee_{\beta}] \subseteq \gee_{\alpha +
    \beta},\ \ \ \alpha,
\beta \in \zz.$
\item[(2)] $[a,b] = -(-1)^{\alpha \beta}[b,a] $ (graded skew-symmetry)
 \item[(3)] $[a,[b,c]] = [[a,b],c] + (-1)^{\alpha \beta}[b,[a,c]]$, \ \ (graded Jacobi
 identity)
\end{description}
for all $a \in \gee_{\alpha}, b \in \gee_{\beta}, c \in \gee_{\gamma}$.
We assume that $\gee$ is nilpotent, so that $\gee^{m}=[\gee,\gee^{m-1}]=0$
for some $m \geq 1$.
The universal enveloping superalgebra $U(\gee)$ is $\zz$-graded and is isomorphic to
$U(\geven) \otimes \bigwedge(\godd)$, where
$\bigwedge(\godd)$ is the exterior (Grassmann) algebra on $\godd$.
We denote
$U(\gee)$ by $\yoo$.

\subsection {Graded-Primitive Ideals}
A \emph{graded-prime} ideal of $\yoo$ is a $\mathbb{Z}_2$-graded
ideal $P$ such that
for any pair $I,J$ of $\zz$-graded ideals of $\yoo$ we have $IJ \subseteq P$
only if $I \subseteq P$ or $J \subseteq P$.
A \emph{graded-primitive} ideal of $\yoo$ is the annihilator of some
simple $\zz$-graded $\yoo$-module.
Let $\textrm{GrSpec}~\yoo$
and
$\textrm{GrPrim}~\yoo$ denote the sets of graded-prime and
graded-primitive
ideals
of $\yoo$, respectively, and let $\textrm{Spec}~U$ and
$\textrm{Prim}~U$ denote the
sets of prime ideals and primitive ideals respectively of the
enveloping algebra $U$ of the Lie algebra $\geven$.
Note that $\textrm{GrPrim}~\yoo \subset \textrm{GrSpec}~\yoo$ and
$\textrm{Prim}~U \subset \textrm{Spec}~U$ ({\cite[3.1.6]{Dixmier:ea96}}).

Below we recount Corollary III of Section 3 in
\cite{Letzter:cmams92}:
\begin{pro}\label{grspec}
Assume that $\gee$ is a finite-dimensional nilpotent Lie superalgebra
over $\mathbb{C}$.
\begin{enumerate}
\item[\rm{(a)}] If $P$ is a graded-prime ideal of $\yoo$, then there exists a
unique
prime ideal $i(P) \in \textrm{Spec}~U$ minimal over $U \cap P$.
\item[\rm{(b)}] The assignment in part (a) produces homeomorphisms of
  topological spaces (relative to the Zariski topology):
\[ i:\textrm{GrSpec}~\yoo \rightarrow \textrm{Spec}~U \]
\[ i:\textrm{GrPrim}~\yoo \rightarrow \textrm{Prim}~U .\]
\end{enumerate}
\end{pro}

The homeomorphism $i$ gives us a bijection between the set of graded-primitive
ideals of $\yoo$ and the primitive ideals of $U$.

\subsection {Induced modules}
Here we relax our assumptions, and
let $\gee$ be an arbitrary Lie superalgebra and $\yoo$ be its
universal enveloping
superalgebra. Let $\hee$ a subsuperalgebra of $\gee$, and $W$ a
$\mathbb{Z}_2$-graded
$\hee$-module, thus a $U(\hee)$-module. The
induced module $Ind_{\hee}^{\gee}W=\yoo \otimes_{U(\hee)}W$ inherits a $\zz$-grading from $W$ and $\yoo$,
and has $\yoo$-action given by
\[ x(u\otimes w) = xu\otimes w,\ x,\ u \in \yoo,\ w \in W \]

Now, we impose the assumption that $\gee$ is a nilpotent Lie superalgebra.

Define
\begin{equation} \Lambda = \{ \lambda \in \gee^* \mid \lambda(\godd)=0\} \end{equation}
\noindent (Then $\Lambda$ can be identified with $\geven^*$).

Two simple $\gee$-modules (resp. $\geven$-modules) are called \emph{weakly
equivalent} if they have
the
same annihilator in $\yoo$ (resp. in $U$).
Let $G_0$ be the group  $\exp(\geven)$.
By
{\cite[Thm.~6.2.4]{Dixmier:ea96}}
there exists a bijective correspondence between the set of
$G_0$-orbits in $\geven^*$ and the set of classes of weakly equivalent
$\geven$-modules.

For $\lambda \in \gee^*$, define $\gee^{\lambda} = \{x \in \gee \mid
\lambda([x,y])=0, \ \forall y \in \gee \}$.
A subsuperalgebra $\kay$ of $\gee$ is said to be \emph{subordinate} to
$\lambda$ if $\lambda([\kay,\kay]) = 0$ and $\gee^{\lambda} \subset \kay$.
A maximal member $\pee$ of the set of subsuperalgebras that are subordinate
to $\lambda$ is called a \emph{polarisation} of $\lambda$.  

Denote by $\{\kay,\lambda\}$ the one-dimensional $\kay$-module given by   
$\lambda$. Thus, $\{\kay,\lambda\}=\mathbb{C}v$ where $x.v=\lambda(x)v$,
for all $x \in \kay$ (i.e. $\{\kay,\lambda\} = \Veven \oplus \Vodd$, where $\Veven = \C v$ and $\Vodd = 0$). Since $\lambda([\kay,\kay])=0$, therefore $\{\kay,\lambda\}$ is a 
well-defined $\kay$-module.
Let $Ind_{\kay}^{\gee}\{\kay,\lambda\}$ denote the induced   
$\gee$-module $\yoo \otimes_{U(\kay)}\{\kay,\lambda\}$.

\bigskip
   
\begin{thm}\label{indmod}{\rm({\cite[Sec. 5.2, Thm. $7'$(b)]{Kac:am77}})}
Let $\gee$ be a finite-dimensional nilpotent Lie superalgebra over an
algebraically closed field of characteristic zero.
Let $\lambda \in \Lambda$,
and let $\pee$ be a polarisation of $\lambda$ with
$\dim \pee_{\bar{0}} = \frac12 \bigl (\dim \geven + \dim(\gee^{\lambda})_{\bar{0}}\bigr )$.
Then the following hold:
\begin{itemize}
\item
The $\gee$-module $\emm_{\lambda} = Ind_{\pee}^{\gee}\{\pee,\lambda\}$
is simple.
\item
The map
$\lambda \rightarrow \emm_{\lambda}$ induces a bijective correspondence
between the set of $G_0$-orbits in $\Lambda$ and
the set of classes of weakly equivalent $\zz$-graded simple $\gee$-modules.
\end{itemize}
\end{thm}

\section{Results on Graded-Primitive Ideals}\label{superres}
Let $\gee$ be a nilpotent Lie superalgebra, $\yoo$ (respectively $U$) be the universal enveloping
superalgebra of $\gee$ (universal enveloping algebra of $\geven$).
Set $\Lambda = \{\lambda \in \gee^* \mid \lambda(\godd) = 0
\}=\geven^*$ as above. An element $\lambda \in \Lambda$ determines the
graded-primitive ideal $P_{\lambda} =
\textrm{ann}_{\yoo}(\mathcal{M}_{\lambda}) \in \textrm{GrPrim}~\yoo$.
Let 
\begin{equation} \Lambda' = \{\lambda \in \Lambda \mid \lambda([\godd,\godd]) = 0\}. \end{equation}
Then we have the following results:
\begin{thm}\label{superres1} Let $\gee = \geven \oplus \godd$ be a nilpotent Lie superalgebra
with universal enveloping superalgebra $\yoo$, and let the universal enveloping algebra of 
$\geven$ be denoted by $U$.
\begin{enumerate}
\item[\rm{(a)}]
Any $\lambda \in \Lambda'$ gives a graded-primitive ideal
$P_{\lambda}$ of $\yoo$ such that $P_{\lambda} \cap U$
is a primitive ideal of $U$. If $\lambda_1,\lambda_2 \in \Lambda'$ are
in the same $G_0$-orbit, then $P_{\lambda_1}=P_{\lambda_2}$.
\item[\rm{(b)}]
If, in addition, we have $[\godd,\godd] = 0$, then
the map $\lambda \rightarrow P_{\lambda}$
induces a bijection between the set of $G_0$-orbits in $\Lambda$ and
the set $\textrm{GrPrim}~\yoo$  of graded-primitive ideals of $\yoo$.
\end{enumerate}
\end{thm}

\subsection{Proof of Theorem~\ref{superres1}}

\textbf{Part \textrm{(a)}}

\noindent \emph{Case 1}:
Suppose $\lambda \in \Lambda'$ is such that
$\lambda([\geven,\geven])=0$. Then we have $\lambda([\gee,\gee])=0$.
So $\gee$ is subordinate to $\lambda$. Thus $\emm_{\lambda} = Ind_{\gee}^{\gee}\{\gee,\lambda\} = \{\gee,\lambda\}$.
Since $\lambda(\godd)=0$, the elements of $\godd$ act trivially on
$\emm_{\lambda}$, so $\emm_{\lambda}$ is a one-dimensional (therefore
simple) $\zz$-graded $\gee$-module (the odd subspace of $\emm_{\lambda}$ is trivial).
Its annihilator in $\yoo$ is the
graded-primitive ideal $P_{\lambda}$ of $\yoo$ generated by all elements
of the form
$x-\lambda(x), \ x \in \gee$. Viewed as a member of $\geven^*$, the
linear map $\lambda$
defines a one-dimensional $\geven$-module $N = \mathbb{C}n$ where
$x.n = \lambda(x)n,$ for all $x \in \geven$.
The primitive ideal of $U$ corresponding to this simple $\geven$-module is
the two sided ideal $Q_{\lambda}$ of $U$ generated by the elements
$\{x-\lambda(x)\vert\
x \in \geven\}$.
It is clear that $Q_{\lambda} = P_{\lambda}\cap U = i(P_{\lambda})$.

\bigskip

\noindent \emph{Case 2}:
Suppose $\lambda \in \Lambda'$ and $\lambda([\geven,\geven]) \neq
0$. Let $\pee \subset \gee$ be a maximal subalgebra subordinate to
$\lambda$. Note that since $\lambda([\godd,\godd])=0$ and
$\lambda([\godd,\geven])
\subset
\lambda(\godd) = 0$, we have $\godd \subset \gee^{\lambda} \subset \pee$.

Thus
$\pee = \pee_{\bar{0}} \oplus \godd$, where $\pee_{\bar{0}}$
 is a maximal subalgebra of
$\geven$ subordinate to $\lambda \ (\in \geven^*)$, which means that
$\dim\pee_{\bar{0}} = \frac12(\dim\geven +
\dim(\gee^{\lambda})_{\bar{0}})$
({\cite[1.12.1]{Dixmier:ea96}}).

Let $\emm_{\lambda} =
Ind_{\pee}^{\gee}\{\pee,\lambda\}$. This is a simple $\gee$-module.
By {\cite[Thm.~6.1.7, Prop.~6.2.3, Thm.~6.2.4]{Dixmier:ea96}} we have
a bijective correspondence
\[ Q_{\lambda} \longleftrightarrow \text{$G_0$-orbit of $\lambda$}\]
between the
primitive ideals of $U$
and the $G_0$-orbits of elements of $\geven^*$.
By Theorem~\ref{indmod} we have a bijective correspondence
\[ \text{$G_0$-orbit of $\lambda$} \longleftrightarrow P_{\lambda} \]
between the set of $G_0$-orbits in $\Lambda$ (= $\geven^*$) and the
set of classes of weakly equivalent $\gee$-modules. Combining the two,
we obtain a correspondence between the primitive ideal $Q_{\lambda}$
of $U$ associated to the $G_0$-orbit of $\lambda$
and the graded-primitive ideal $P_{\lambda}$ of $\yoo$ that is
the annihilator in $\yoo$ of all the
simple $\gee$-modules weakly equivalent to $\emm_{\lambda}$.
Note that $Q_{\lambda} = \textrm{ann}_U(\mathcal{N}_{\lambda})$, where
$\mathcal{N}_{\lambda} =
Ind_{\pee_{\bar{0}}}^{\geven}\{\pee_{\bar{0}},\lambda\}$, which is a
simple $\geven$-module (from {\cite[Thm. 6.1.1]{Dixmier:ea96}}).

Now, $\pee = \pee_{\bar{0}} \oplus \godd$,
so we may choose linearly independent elements $e_1,e_2,\dots,e_r \in \geven$ such that $\geven = \langle
e_1,e_2,\dots,e_r \rangle \oplus \pee_{\bar{0}}$ and $\gee = \langle e_1,e_2,\dots,e_r \rangle \oplus
\pee$.

Let $\{\pee,\lambda\} = \mathbb{C}v$ and $\{\pee_{\bar{0}},\lambda\} =
\mathbb{C}w$.
Then, the $\gee$-module $\emm_{\lambda}$ consists of linear combinations of
elements of the type $e_1^{a_1}e_2^{a_2}\cdots e_r^{a_r}\otimes v$ and
the $\geven$-module $\mathcal{N}_{\lambda}$ consists of linear combinations of elements
of the type $e_1^{a_1}e_2^{a_2}\cdots e_r^{a_r}\otimes w$, where $a_1,a_2,
\dots,a_r$ are non-negative integers.

Since $\godd$ acts trivially on $\emm_{\lambda}$, we have that
$\emm_{\lambda}$ viewed as a
$\geven$-module is annihilated by the same elements in $U$ as $\mathcal{N}_{\lambda}$.
Therefore $Q_{\lambda} = P_{\lambda}\cap U = i(P_{\lambda})$.
Let $S_{\lambda}$ denote a minimal set of generators of the primitive ideal $Q_{\lambda}$
of $U$ (we can do this because $U$ is Noetherian),
 and let $T_{\lambda}$ denote a minimal set of generators of the graded-primitive
ideal $P_{\lambda}$ of $\yoo$. Since $\godd$ acts trivially on $\emm_{\lambda}$,
we may assume $T_{\lambda}$ contains the basis elements $\{f_1,\dots,f_s\}$, of a
fixed basis of $\godd$.
In fact, we can choose $T_{\lambda}$ to be such that $T_{\lambda} = S_{\lambda} \cup \{f_1,\dots,f_s\}$.

If $\lambda_1, \lambda_2 \in \Lambda'$ are in the same $G_0$-orbit,
then $Q_{\lambda_1}=Q_{\lambda_2}$. Then we may suppose $S$ is
a minimal set of generators of the
primitive ideal $Q_{\lambda_1}=Q_{\lambda_2}$ of $U$. Then $T_{\lambda_1} = S \cup \{f_1,\dots,f_s\}=T_{\lambda_2}$ is a minimal set of  generators of
the graded-primitive ideal $P_{\lambda_i}$ of $\yoo$.
Therefore
$P_{\lambda_1} = P_{\lambda_2}$.

This completes the proof of Part \textrm{(a)} of Theorem~\ref{superres1}.
\bigskip

\textbf{Part \textrm{(b)}}

If $[\godd,\godd] = 0$, then $\Lambda' = \Lambda$.
By {\cite[Sec.~6.1.5, Thm.~6.1.7]{Dixmier:ea96}},
each primitive ideal
of $U$ is
of the form $Q_{\lambda}$, where $\lambda \in \geven^*$. Moreover,
if $\kay$ is any
maximal subalgebra of $\geven$ subordinate to $\lambda$, then
$Q_{\lambda} = \textrm{ann}_U(Ind_{\kay}^{\geven}\{\kay,\lambda\})$.
Thus, the ideal $Q_{\lambda}$ depends only on $\lambda$.

The map $\lambda \rightarrow Q_{\lambda}$ gives a bijection between the set of
$G_0$-orbits in $\geven^*$ and $\textrm{Prim}~U$ ({\cite[Thm.~6.2.4]{Dixmier:ea96}}).

We see from Cases 1 and 2 above that any $\lambda \in \Lambda'=\Lambda$ determines
a graded-primitive ideal $P_{\lambda}$ of $\yoo$ such that $i(P_{\lambda})
= P_{\lambda} \cap U \in \textrm{Prim}\;U$.
Here, $\Lambda' = \Lambda$, so we can replace $\Lambda'$ by $\Lambda$
in the statement above.
So, from Proposition \ref{grspec} and Theorem \ref{indmod}, the map
$\lambda \rightarrow P_{\lambda}$ induces a bijection between the set of
graded-primitive ideals of $\yoo$ and the set of $G_0$-orbits in $\Lambda$.
This completes the proof of Part \textrm{(b)} of Theorem \ref{superres1}.  \qed

\bigskip

Before going on, we recall some well-known
results:
\bigskip

From {\cite[Thm. 4.7.9, Sec 4.7.10, and Prop. 6.2.2]{Dixmier:ea96}},
we know that for $\lambda \in \geven^*$ that $U/Q_{\lambda} \simeq
\mathcal{A}_n$, the $n$-th Weyl algebra (Definition~\ref{D:weylalg}), where
$2n$ is the rank of the bilinear form $B_{\lambda}$ on $\geven$ defined by
$B_{\lambda}(x,y)=\lambda([x,y])$. Thus there exist elements
$x_i,y_i,\ 1\leq i\leq n$, in $\geven$, such that $X_i = x_i +
Q_{\lambda}\ Y_i = y_i + Q_{\lambda}$ satisfy the Weyl
relations
in $U/Q_{\lambda}$:
\[ [X_i,X_j] = 0 = [Y_i,Y_j] \]
\[ [X_i,Y_j] = \lambda([X_i,Y_j]) = \delta_{ij}1. \]

The number $n$ is called the \emph{weight} of the
primitive ideal
$Q_{\lambda}$.
 What can we say about the factor $\yoo /P_{\lambda}$? Let us recall the
 following result from \cite{BellMusson:lms90}:

\begin{thm}\label{primfac}{\rm({\cite[Cor. B]{BellMusson:lms90}})}
Suppose $k$ is an algebraically closed field of characteristic
zero, $\gee$ is a finite-dimensional nilpotent Lie superalgebra over
$k$, and $\yoo = U(\gee)$ is the universal enveloping superalgebra of $\gee$.
If $P$ is a primitive ideal of
$\yoo$, then
\[ \yoo /P \simeq M_s(\mathcal{A}_n(k)); \] and if $P$ is a graded-primitive ideal of $\yoo$, then
\[ \yoo /P \simeq M_s(\mathcal{A}_n(k)) \ \text{or} \ \yoo /P \simeq
 M_s(\mathcal{A}_n(k))\times M_s(\mathcal{A}_n(k)), \]
where $s = 2^m$, $m,n$ are non-negative integers, and $M_s(\mathcal{A}_n(k))$
denotes the algebra of $s \times s$ matrices over the $n$-th Weyl
algebra $\mathcal{A}_n(k)$.
\end{thm}

Therefore, as a consequence of Theorem~\ref{superres1}, we have the following result:
\begin{thm}\label{superres2} Let $\gee$ be a finite-dimensional nilpotent Lie superalgebra over
$\C$ with universal enveloping superalgebra $\yoo$.
\begin{enumerate}
\item[\rm{(a)}]
Let $\lambda \in \Lambda'$ and suppose
$P_{\lambda}$ is the corresponding graded-primitive ideal of
$\yoo$. Then,
\[ \yoo/P_{\lambda} \simeq \mathcal{A}_n \]
\noindent where $2n=\text{rank}(B_{\lambda})$ on $\geven$, and
$\mathcal{A}_n$ is the $n$-th Weyl algebra over
$\mathbb{C}$. 
\item[\rm{(b)}]
Suppose $\gee = \geven \oplus \godd$ satisfies the condition
$[\godd,\godd]=0$.
Then, for any $P \in \textrm{GrPrim}~\yoo$, we have
\[ \yoo/P \simeq \mathcal{A}_n \]
\noindent for a unique non-negative integer $n$.
\end{enumerate}
\end{thm}
\subsection{Proof of Theorem~\ref{superres2}}

\noindent \textbf{(a)} Suppose $P_{\lambda}$ is the graded primitive ideal corresponding to
$\lambda \in \Lambda'$.
Then the factor $\yoo /P_{\lambda}$ is constructed by factoring $\yoo$ by
the relations $T_{\lambda} = 0$, where $T_{\lambda}$ is a minimal set of generators of
$P_{\lambda}$. But we can assume that $T_{\lambda} = S_{\lambda} \cup \{f_1,\dots,f_s\}$, where
$S_{\lambda}$ is a minimal set of generators of the primitive ideal
$Q_{\lambda}$ of $U$, and $\{f_1,\dots,f_s\}$ is a fixed basis of $\godd$.
So, taking the relevant Poincar\'{e}-Birkhoff-Witt basis of $\yoo$, we see that
$\yoo /P_{\lambda} \simeq U/Q_{\lambda} \simeq \mathcal{A}_n$.
In the notation of Theorem~\ref{primfac}, we have
$s = 1$ and $n$ is the weight of the primitive ideal of $U$ that is
in one-one correspondence with $P \in \textrm{GrPrim}~\yoo$.

\bigskip

\noindent \textbf{(b)}
If $[\godd,\godd] = 0$, then, by Part (\textrm{b}) of Theorem~\ref{superres1}, any
$P \in \textrm{GrPrim}~\yoo$ is of the form $P_{\lambda}$ for some $\lambda \in \Lambda'$.
Therefore, by Part (a) above, we have $\yoo / P =\yoo /
P_{\lambda}\simeq 
\mathcal{A}_n$, where $2n = \text{rank}\ B_{\lambda}$ on $\geven$.
This proves Part (b).  \qed

\section{Applications} 
\subsection{$\gee = \gl(m,n)^+, m \neq n$}
Let $\gee = \gl(m,n)^+$, where $m \geq 1,\ n \geq 1$, and $m \neq n$.
Then $\gee$ is the nilpotent Lie superalgebra of strictly upper
triangular matrices in the general linear Lie
superalgebra $\gl(m,n)$ over $\mathbb{C}$.
 In matrix notation, $\gee$ is defined to
be the set
of block matrices
\[
  \begin{pmatrix}
   A & B\\
   0 & D
  \end{pmatrix} ,
\]

\noindent where $A,D$ are strictly upper triangular matrices of sizes
$m \times m$
and $n \times n$, respectively, and $B$ is an arbitrary
$m \times n$ matrix.

The Lie superbracket on $\gee$ is as follows:
\[
\left[
  \begin{pmatrix}
   A & B\\
   0 & D
  \end{pmatrix},
  \begin{pmatrix}
   A' & B'\\
   0 &  D'
  \end{pmatrix}
\right]
=
\begin{pmatrix}
AA'-A'A & BD'-B'D+AB'-A'B\\
0 &   DD'-D'D
\end{pmatrix}
\]
  
The even part
\[
   \geven =
   \left\{
    \begin{pmatrix}
     A & 0\\
     0 & D
    \end{pmatrix} \mid A \in \en_m,\ D \in \en_n
    \right\}
\]

\noindent is isomorphic to the nilpotent Lie algebra $\en_m \times
\en_n$, while
the odd part
\[
 \godd =  
\left\{ 
\begin{pmatrix} 
 0 & B\\
 0& 0
\end{pmatrix} \mid B \in {M}_{m \times n}(\C)
\right\}
\]
 
\noindent has dimension $mn$ and $[\godd,\godd]=0$. (Note that $\en_1
= 0$.)

\begin{defn}\label{range}
We define an integer-valued function $s_i$ for $i = m$ or
$n$, as follows:
\[ s_i= \begin{cases}
     \frac14(i-2)i & \text{if $i$ is even},\\
     \frac14(i-1)^2 & \text{if $i$ is odd}.
  \end{cases} \]
  \end{defn}
     
Again, by Theorem~\ref{superres1}, there is a bijection between the set
of $G_0$-orbits in $\Lambda$ and the set $\textrm{GrPrim}~\yoo$
 of graded-primitive ideals of
$\yoo$. By Theorem~\ref{superres2}, for any $P \in
\textrm{GrPrim}~\yoo$,
the quotient $\yoo/P \simeq A_r$ where $r$ is the weight
of the ideal $Q = P \cap U \in \textrm{Prim}~U$.
In this case, $U$ is the enveloping algebra of the nilpotent Lie
algebra $\geven = \en_m \times \en_n$. Let $U_m$ and $U_n$ be the
universal enveloping
algebras of the Lie algebras $\en_m$ and $\en_n$, respectively.
By Corollary 4.10 in {\cite{Mukherjee:thesis04}}, the weights of members of
$\textrm{Prim}~U_m$ range through
$0,1,\dots,s_m$, and the weights of members of
$\textrm{Prim}~U_n$ range through $0,1,\dots,s_n$.

For $\lambda \in \Lambda$, any subsuperalgebra $\pee$ of $\gee$ maximally
subordinate to $\lambda$ is of the form
$\pee_{\bar{0}} \oplus \godd$, where $\pee_{\bar{0}}$ is a
polarisation in $\geven$ of $\lambda \ (\in \geven^*)$, as described
in the proof of Theorem~\ref{superres1}. The subalgebra
$\pee_{\bar{0}}$
can be chosen to be of the form $\hee_m \times \hee_n$, where
$\hee_m$ is a polarisation of $\lambda \vert_{\en_{m}}$ and
$\hee_n$ is a polarisation of $\lambda \vert_{\en_{n}}$.
So the codimension of $\pee_{\bar{0}}$ in $\geven$ is $r_m + r_n$,
where $r_m$ can range through $0,1,\dots,s_m$ and $r_n$
can range through $0,1,\dots,s_n$ (see remark at the end of Chapter 2 in
{\cite{Mukherjee:thesis04}}).
     
Thus, we have the following:
\begin{pro}
If $\yoo$ is the universal enveloping superalgebra of
$\gee = \gl(m,n)^+$,
then, for any $P \in \textrm{GrPrim}~\yoo$, the quotient
$\yoo /P \simeq \mathcal{A}_{r_m+r_n}$, where $r_m$ and $r_n$ are
unique non-negative integers; $0 \leq r_m \leq s_m$, and
$0 \leq r_n \leq s_n$, and $s_{m},\ s_{n}$ are given by Definition~\ref{range}
\end{pro}
\subsection{The Heisenberg Lie Superalgebra}
The next example, which comes from {\cite[Sec. 0.2(a)]{BellMusson:lms90}}, shows
that Theorem~\ref{superres2} may not hold when
$[\godd,\godd] \neq 0$.

Let $\gee$ be the nilpotent Lie superalgebra over $\mathbb{C}$
with basis for $\geven$ given by
$x,y,z$ and basis for $\godd$ given by $a,b$. Let all Lie
superbrackets be zero except $[x,y]=z=-[y,x]$ and $[a,b]=z=[b,a]$. Thus,
$[\geven,\geven]=[\godd,\godd]=\mathbb{C}z$.
Let $\lambda \in \gee^*$ be such that $\lambda(\godd)=0$ and
$\lambda(z)=1$.
Then $\gee^{\lambda}=\mathbb{C}z$. The basis elements $x,a,z$ span a
subsuperalgebra $\pee$ that is subordinate to $\lambda$ and is of maximal
dimension. Let $\{\pee,\lambda\}=\mathbb{C}v$ denote the   
one-dimensional $\pee$-module given by $\lambda$.
By Theorem~\ref{indmod}, the $\mathbb{Z}_2$-graded $\gee$-module  
$\emm_{\lambda} = Ind_{\pee}^{\gee}\{\pee,\lambda\}$ is irreducible.
We can see that $\emm_{\lambda}$ is spanned by elements of the form
\[ y^r b^s\otimes v \]
\noindent where $r,s$ are non-negative integers, $r \geq 0$ and
$s=0$ or $1$. The annihilator in $\yoo$ of this module is the
graded-primitive ideal $P_{\lambda}$ generated by $z-1$.
But we see that $\yoo/P_{\lambda} \simeq
M_2(\mathbb{C})\otimes_{\mathbb{C}}\mathcal{A}_1 \simeq
M_2(\mathcal{A}_1)$,
because in the quotient
$\yoo/P_{\lambda}$, the elements $\bar{x}$ and $\bar{y}$
generate a copy of $\mathcal{A}_1$ and the elements
$\bar{1}, \bar{a}, \bar{b}, \bar{a}\bar{b}$ form a basis for
the $\mathbb{C}$-algebra $M_2(\mathbb{C})$. 

\section{Conclusion}

In this work we have used results of Benoist, Fernando, Kirillov and Dixmier  to study
modules and coadjoint orbits of finite-dimensional nilpotent Lie algebras. We have also derived
results about
certain kinds of simple infinite-dimensional modules and the corresponding graded-primitive
ideals of the universal enveloping superalgebra of nilpotent Lie superalgebra, using the work of Bell, Musson, Letzter, and
Kac. Our investigations suggest the following problems for
future study.

\begin{enumerate}
\item Let $\gee$ be a finite-dimensional nilpotent Lie algebra,
let $f \in \gee^*$,
let $\Omega_f$ be the coadjoint orbit containing $f$,
 and let $\pee \in P(f)$ be a polarisation. If the dimension of $\Omega_f$ is greater
than two, when  is the variety $\Omega_f \cap (f + \pee^{\top})$
lagrangian? How many elements does the set
$\mathcal{E}_I^{\pee,f}$ of simple modules have?

\item
Let $\gee$ be a finite-dimensional nilpotent Lie algebra over $\C$
with an involution $\sigma$ (an automorphism of order 2) , and let $\hee$ be the set of fixed points of
$\sigma$. A simple $\gee$-module is said to be $\sigma$-\emph{spherical} if
it contains a nontrivial vector annihilated by $\hee$.
Let $j$ be the principal anti-automorphism of the universal enveloping
algebra $U$ of $\gee$ such that $x \rightarrow -x$ for all $x \in \gee$.
Let $\textrm{Prim}~U$ be the set of primitive ideals of $U$. Set
$\textrm{Prim}_{\sigma}U = \{I\in \textrm{Prim}~U \mid I^{\sigma}
=I^{j}\}$. In {\cite{Benoist:cm90}}, Benoist showed for $\sigma$, a fixed involution on $\gee$, that
 there is a bijection
between the set of isomorphism classes of $\sigma$-spherical simple
$\gee$-modules and the
set of ideals $\textrm{Prim}_{\sigma}U$,
and he also gave a classification and several
constructions of these modules.

Is it possible to develop a theory of
a $\sigma$-spherical simple modules for finite-dimensional nilpotent Lie
\emph{superalgebras}
 $\gee$ with an involution?
Are these modules  in one-to-one correspondence with a subset of
$\textrm{GrPrim}~\yoo$, the set of graded-primitive ideals of the universal
enveloping superalgebra of $\gee$ ?
\end{enumerate}

 \end{document}